\documentclass[10pt]{article}

\usepackage{scrpage2}
\clearscrheadfoot
\usepackage{amsmath}
\usepackage{amsfonts}
\usepackage{amssymb}

\usepackage{xcolor}
\usepackage{graphicx}
\usepackage{tikz}
\usepackage{caption}
\usepackage{subcaption}

\usepackage{hyperref}
\hypersetup{linktocpage,colorlinks=true,pdfborder={0 0 0},linkcolor=red,urlcolor=red}

\usepackage{geometry}
\title{Minimal completely asymmetric $(4;n)$-regular matchstick graphs}
\author{\\Mike Winkler\footnote{Fakult\"at f\"ur Mathematik, Ruhr-Universit\"at Bochum, Germany, mike.winkler@ruhr-uni-bochum.de}\quad Peter Dinkelacker\footnote{Togostr. 79, 13351 Berlin, Germany, peter@grity.de}\quad Stefan Vogel\footnote{Raun, Dorfstr. 7, 08648 Bad Brambach, Germany, backebackekuchen16@gmail.com}\\ \\}
\date{May 1, 2018}

\begin{document}
  
  \maketitle
  
  \begin{abstract}
	\noindent A matchstick graph is a graph drawn with straight edges in the plane such that the edges have unit length, and non-adjacent edges do not intersect. We call a matchstick graph $(m;n)$-regular if every vertex has only degree $m$ or $n$. In this article we present the latest known $(4;n)$-regular matchstick graphs for $4\leq n\leq11$ with a minimum number of vertices and a completely asymmetric structure.
	\\ \\
    We call a matchstick graph \textit{completely asymmetric}, if the following conditions are complied.
    
    \begin{itemize}
      \item The graph is rigid.
      \item The graph has no point, rotational or mirror symmetry.
      \item The graph has an asymmetric outer shape.
      \item The graph can not be decomposed into rigid subgraphs and rearrange to a similar graph which contradicts to any of the other conditions.
    \end{itemize}
  
  \end{abstract}
  
  \newpage
  
  \section{\large{Introduction}}
  
  A matchstick graph is a planar unit-distance graph. That is a graph drawn with straight edges in the plane such that the edges have unit length, and non-adjacent edges do not intersect. We call a matchstick graph $(m;n)$-regular if every vertex has only degree $m$ or $n$.
  \\ \\
  For $m\leq n$ minimal $(4;n)$-regular matchstick graphs with a minimum number of vertices only exist for $4\leq n\leq11$. The smallest known $(4;n)$-regular matchstick graph for $n=4$, also named 4-regular, is the so called \textit{Harborth graph} consisting of 52 vertices and 104 edges. Except for $n=10$ and $n=11$ all currently smallest known $(4;n)$-regular matchstick graphs for $4\leq n\leq11$ are symmetric. They consist of $104, 115, 117, 159, 126, 273, 231$ and $771$ edges and were presented in an earlier paper by the authors [6].
  \\ \\
  It is an open problem how many different $(4;n)$-regular matchstick graphs with a minimum number of vertices for $4\leq n\leq11$ exist and which is the least minimal number. "Our knowledge on matchstick graphs is still very limited. It seems  to be hard to obtain rigid mathematical results about them. Matchstick problems constructing the minimal example can be quite challenging. But the really hard task is to rigidly prove that no smaller example can exist." [2]
  \\ \\
  In this article we present the $(4;n)$-regular matchstick graphs for $4\leq n\leq11$ with the smallest currently known number of vertices and a completely asymmetric structure. The definition of a \textit{completely asymmetric matchstick graph} is given in Chapter 3. The graphs were discovered in the days from March 17, 2016 -- April 15, 2018 and were presented for the first time in a German mathematics internet forum [3]. The graphs for $n=10$ and $n=11$ are also the currently smallest known examples.
  \\ \\
  The geometry, rigidity or flexibility of the graphs in this article has been verified by Stefan Vogel with a computer algebra system named \textsc{Matchstick Graphs Calculator} (MGC) [5]. This remarkable software created by Vogel runs directly in web browsers. A special version of the MGC contains all graphs from this article and is available under this \href{http://mikewinkler.co.nf/matchstick_graphs_calculator.htm}{weblink}\footnote{http://mikewinkler.co.nf/matchstick\_graphs\_calculator.htm}. The method Vogel used for the calculations he describes in a separate German article [4].
  \\ \\
  \textbf{Remark}: \textit{Proofs for the existence of the graphs shown in this article}.
  
  \begin{enumerate}
    \item The MGC contains a constructive proof for each graph. We are using this online reference, because these proofs are too extensive to reproduce here.
  \end{enumerate}
  \textbf{Note}: In the PDF version of this article the vector graphics can be viewed with the highest zoom factor to see the smallest details. For example the very small rhombus in Figure 13.
  
  \newpage
  
  \section{\large{Rigid subgraphs}}
  
  The geometry of the graphs in this paper is based on two types of rigid subgraphs we call the \textit{kite} (Fig. 1 a) and the \textit{triplet kite} (Fig. 1 b). The kite is a $(2;4)$-regular matchstick graph consisting of 12 vertices and 21 edges and has a vertical symmetry. The triplet kite is a $(2;3;4)$-regular matchstick graph consisting of 22 vertices and 41 edges and has a vertical symmetry.
  \\ \\
  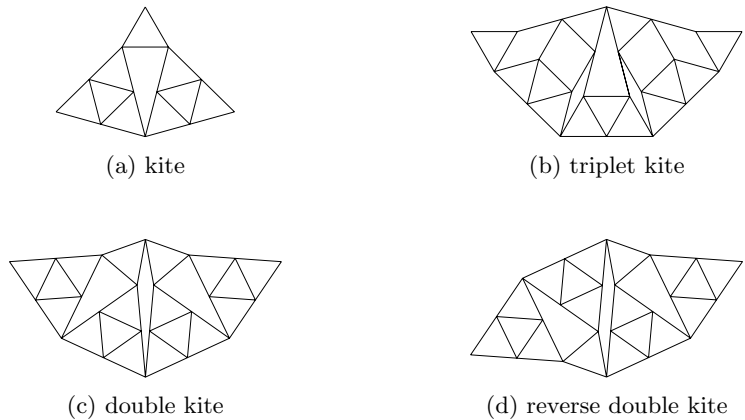
\begin{figure}[!ht]
    \centering
    \begin{minipage}[t]{0.45\linewidth}
      \centering
      \begin{tikzpicture}
      [y=0.4pt, x=0.4pt, yscale=-1.0, xscale=1.0]
      \draw[line width=0.01pt]
      (86.5346,1050.0456) -- (128.7085,1038.3319)
      (139.6510,995.9514) -- (128.7085,1038.3319)
      (128.7085,1038.3319) -- (97.4772,1007.6651)
      (97.4772,1007.6651) -- (86.5346,1050.0456)
      (97.4772,1007.6651) -- (139.6510,995.9514)
      (139.6510,995.9514) -- (108.4198,965.2847)
      (108.4198,965.2847) -- (97.4772,1007.6651)
      (108.4198,965.2847) -- (64.6494,965.2847)
      (86.5346,927.3784) -- (64.6494,965.2847)
      (64.6494,965.2847) -- (86.5346,1050.0456)
      (86.5346,1050.0456) -- (2.1869,1026.6183)
      (2.1869,1026.6183) -- (64.6494,965.2847)
      (75.5920,1007.6651) -- (33.4182,995.9515)
      (33.4182,995.9515) -- (44.3608,1038.3320)
      (44.3608,1038.3320) -- (75.5920,1007.6651)
      (128.7085,1038.3319) -- (170.8824,1026.6183)
      (170.8824,1026.6183) -- (139.6510,995.9514)
      (108.4198,965.2847) -- (86.5346,927.3784);
      \end{tikzpicture}
      \subcaption{kite}
    \end{minipage}
    \quad
    \begin{minipage}[t]{0.45\linewidth}
      \centering
      \begin{tikzpicture}
      [y=0.4pt, x=0.4pt, yscale=-1.0, xscale=1.0]
      \draw[line width=0.01pt]
      (119.7716,968.4091) -- (97.8865,1006.3153)
      (174.4846,1048.6958) -- (86.9439,1048.6958)
      (108.8290,1010.7896) -- (152.5994,1010.7896)
      (130.7142,1048.6958) -- (108.8290,1010.7896)
      (119.7716,968.4091) -- (108.8290,1010.7896)
      (86.9439,1048.6958) -- (108.8290,1010.7896)
      (130.7142,1048.6958) -- (152.5994,1010.7896)
      (152.5994,1010.7896) -- (174.4846,1048.6958)
      (97.8865,1006.3153) -- (86.9439,1048.6958)
      (24.4813,987.3622) -- (86.9439,1048.6958)
      (97.8865,1006.3153) -- (66.6552,975.6486)
      (55.7126,1018.0290) -- (97.8865,1006.3153)
      (55.7126,1018.0290) -- (66.6552,975.6486)
      (66.6552,975.6486) -- (24.4813,987.3622)
      (88.5404,937.7423) -- (66.6552,975.6486)
      (2.5962,949.4560) -- (46.3665,949.4560)
      (46.3665,949.4560) -- (24.4813,987.3622)
      (88.5404,937.7423) -- (46.3665,949.4560)
      (88.5404,937.7423) -- (119.7716,968.4091)
      (141.6568,968.4091) -- (163.5420,1006.3153)
      (141.6568,968.4091) -- (152.5994,1010.7896)
      (163.5420,1006.3153) -- (174.4846,1048.6958)
      (141.6568,968.4091) -- (152.5994,1010.7896)
      (236.9471,987.3622) -- (174.4846,1048.6958)
      (163.5420,1006.3153) -- (194.7732,975.6485)
      (205.7158,1018.0290) -- (163.5420,1006.3153)
      (205.7158,1018.0290) -- (194.7732,975.6485)
      (194.7732,975.6485) -- (236.9471,987.3622)
      (172.8881,937.7423) -- (194.7732,975.6485)
      (215.0619,949.4560) -- (236.9471,987.3622)
      (172.8881,937.7423) -- (215.0619,949.4560)
      (141.6568,968.4091) -- (172.8881,937.7423)
      (119.7716,968.4091) -- (130.7142,926.0286)
      (130.7142,926.0286) -- (141.6568,968.4091)
      (2.5962,949.4560) -- (24.4813,987.3622)
      (215.0619,949.4560) -- (258.8323,949.4560)
      (258.8323,949.4560) -- (236.9471,987.3622)
      (88.5404,937.7423) -- (130.7142,926.0286)
      (130.7142,926.0286) -- (172.8881,937.7423);
      \end{tikzpicture}
      \subcaption{triplet kite}
    \end{minipage}
    
    \begin{center}\end{center}
    
    \begin{minipage}[t]{0.45\linewidth}
      \centering
      \begin{tikzpicture}
      [y=0.4pt, x=0.4pt, yscale=-1.0, xscale=1.0]
      \draw[line width=0.01pt]
      (54.8439,1011.4230) -- (30.2950,975.1850)
      (49.4036,935.8060) -- (30.2950,975.1850)
      (30.2950,975.1850) -- (73.9525,972.0441)
      (73.9525,972.0441) -- (54.8439,1011.4230)
      (73.9525,972.0441) -- (49.4036,935.8060)
      (49.4036,935.8060) -- (93.0611,932.6651)
      (93.0611,932.6651) -- (73.9525,972.0441)
      (93.0611,932.6651) -- (126.4125,961.0117)
      (134.2857,917.9552) -- (126.4125,961.0117)
      (126.4125,961.0117) -- (54.8439,1011.4230)
      (54.8439,1011.4230) -- (134.2857,1048.1977)
      (134.2857,1048.1977) -- (126.4125,961.0117)
      (90.6282,986.2174) -- (130.3491,1004.6047)
      (130.3491,1004.6047) -- (94.5648,1029.8103)
      (94.5648,1029.8103) -- (90.6282,986.2174)
      (30.2950,975.1850) -- (5.7461,938.9469)
      (5.7461,938.9469) -- (49.4036,935.8060)
      (93.0611,932.6651) -- (134.2857,917.9552)
      (213.7275,1011.4230) -- (238.2764,975.1850)
      (219.1678,935.8060) -- (238.2764,975.1850)
      (238.2764,975.1850) -- (194.6189,972.0441)
      (194.6189,972.0441) -- (213.7275,1011.4230)
      (194.6189,972.0441) -- (219.1677,935.8060)
      (219.1678,935.8060) -- (175.5102,932.6651)
      (175.5102,932.6651) -- (194.6189,972.0441)
      (175.5102,932.6651) -- (142.1589,961.0117)
      (134.2857,917.9552) -- (142.1589,961.0117)
      (142.1589,961.0117) -- (213.7275,1011.4230)
      (213.7275,1011.4230) -- (134.2857,1048.1977)
      (134.2857,1048.1977) -- (142.1589,961.0117)
      (177.9432,986.2174) -- (138.2223,1004.6047)
      (138.2223,1004.6047) -- (174.0066,1029.8103)
      (174.0066,1029.8103) -- (177.9432,986.2174)
      (238.2764,975.1850) -- (262.8253,938.9469)
      (262.8253,938.9469) -- (219.1677,935.8060)
      (175.5102,932.6651) -- (134.2857,917.9552);
      \end{tikzpicture}
      \subcaption{double kite}
    \end{minipage}
    \quad
    \begin{minipage}[t]{0.45\linewidth}
      \centering
      \begin{tikzpicture}
      [y=0.4pt, x=0.4pt, yscale=-1.0, xscale=1.0]
      \draw[line width=0.01pt]
      (53.0582,956.1585) -- (92.7791,937.7712)
      (128.5634,962.9769) -- (92.7791,937.7712)
      (92.7791,937.7712) -- (88.8425,981.3641)
      (88.8425,981.3641) -- (53.0582,956.1585)
      (88.8425,981.3641) -- (128.5634,962.9769)
      (128.5634,962.9769) -- (124.6268,1006.5698)
      (124.6268,1006.5698) -- (88.8425,981.3641)
      (124.6268,1006.5698) -- (91.2754,1034.9164)
      (132.5000,1049.6262) -- (91.2754,1034.9164)
      (91.2754,1034.9164) -- (53.0582,956.1585)
      (53.0582,956.1585) -- (3.9604,1028.6346)
      (3.9604,1028.6346) -- (91.2754,1034.9164)
      (72.1668,995.5375) -- (47.6179,1031.7755)
      (47.6179,1031.7755) -- (28.5093,992.3965)
      (28.5093,992.3965) -- (72.1668,995.5375)
      (92.7791,937.7712) -- (132.5000,919.3839)
      (176.1575,987.6460) -- (211.9418,1012.8517)
      (176.1575,987.6460) -- (136.4366,1006.0333)
      (136.4366,1006.0333) -- (140.3732,962.4403)
      (140.3732,962.4403) -- (176.1575,987.6460)
      (140.3732,962.4403) -- (173.7245,934.0937)
      (140.3732,962.4403) -- (132.5000,919.3839)
      (132.5000,919.3839) -- (173.7245,934.0937)
      (173.7245,934.0937) -- (211.9418,1012.8517)
      (211.9418,1012.8517) -- (261.0396,940.3756)
      (261.0396,940.3756) -- (173.7246,934.0937)
      (192.8332,973.4727) -- (217.3821,937.2346)
      (217.3821,937.2346) -- (236.4907,976.6136)
      (236.4907,976.6136) -- (192.8332,973.4727)
      (211.9418,1012.8517) -- (172.2209,1031.2390)
      (172.2209,1031.2390) -- (176.1575,987.6460)
      (172.2209,1031.2390) -- (136.4366,1006.0333)
      (124.6268,1006.5698) -- (132.5000,1049.6262)
      (172.2209,1031.2390) -- (132.5000,1049.6262)
      (132.5000,919.3839) -- (128.5634,962.9768)
      (136.4366,1006.0333) -- (132.5000,1049.6262);
      \end{tikzpicture}
      \subcaption{reverse double kite}
    \end{minipage}
    
    \begin{center}\end{center}
    
    \caption{Rigid subgraphs}
  \end{figure}
  
  \noindent
  \\
  Two kites can be connected to each other in two useful ways. We call these subgraphs the \textit{double kite} (Fig. 1 c) and the \textit{reverse double kite} (Fig. 1 d), both consisting of 22 vertices and 42 edges. What makes the subgraphs (c) and (d) so useful is the fact that they have only two vertices of degree 2. Two of these subgraphs can be used to connect two vertices of degree 2 at different distances by using them like clasps. This property has been used for $n=8,\dots,11$ (Fig. 9 -- 12).
  \\ \\
  The geometry of the graphs for $n=4,\dots,7$ is based on the triplet kite. The geometry of the graphs for $n=8,\dots,10$ is based on the kite. The geometry of the graph for $n=11$ is a special one, but also this graph contains kites.
  
  \newpage
  
  \section{Completely asymmetric matchstick graphs}
  
  We call a matchstick graph \textit{completely asymmetric}, if the following conditions are complied.
  \\ \\
  $\bullet$ The graph is rigid.\\
  $\bullet$ The graph has no point, rotational or mirror symmetry.\\
  $\bullet$ The graph has an asymmetric outer shape.\\
  $\bullet$ The graph can not be decomposed into rigid subgraphs and rearrange to a similar graph which contradicts to any of the other conditions.
  \\ \\
  \noindent For example the next two matchstick graphs are asymmetric, but not \textit{completely} asymmetric, because only three of the four conditions are complied.
  
  \begin{figure}[!ht]
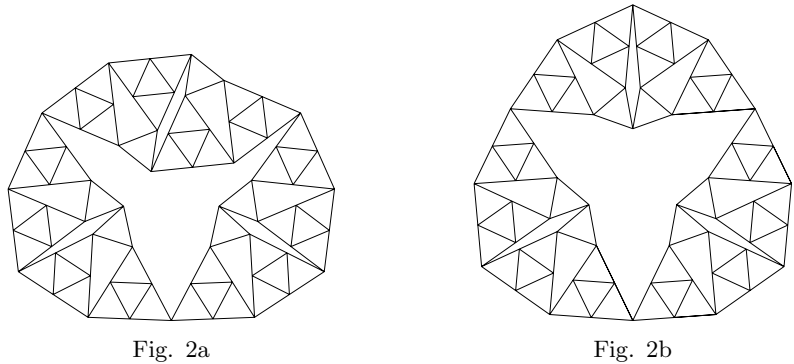

    \centering
    \begin{minipage}[t]{0.45\linewidth}
     \centering

     \subcaption*{Fig. 2a}
   \end{minipage}
   \quad
   \begin{minipage}[t]{0.45\linewidth}
    \centering
     %
     \subcaption*{Fig. 2b}
    \end{minipage}
    \caption{4-regular matchstick graphs with 63 vertices and 126 edges. Figure 2a shows an asymmetric graph consisting only of one type of rigid subgraph - the kite. And six kites can be rearrange to a symmetric graph with a rotational symmetry of order 3 (Fig. 2b). Therefore the graph in Figure 2a is asymmetric, but not completely asymmetric.}
    \end{figure}
  
    \begin{figure}[!ht]
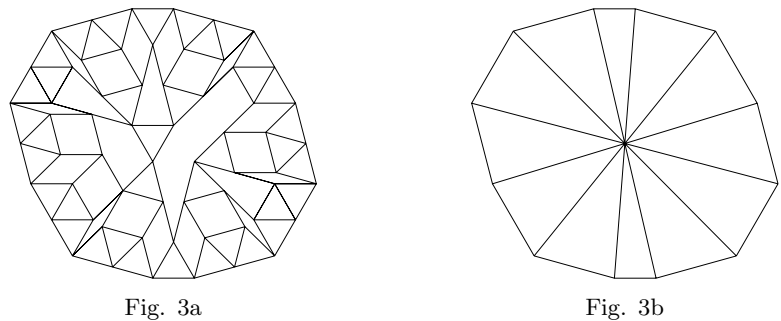

    \centering
    \begin{minipage}[t]{0.45\linewidth}
     \centering
     %
      \subcaption*{Fig. 3a}
     \end{minipage}
     \quad
     \begin{minipage}[t]{0.45\linewidth}
      \centering
     %
    \subcaption*{Fig. 3b}
    \end{minipage}
    \caption{Figure 3a shows a $(4;5)$-regular matchstick graph with 60 vertices and 121 edges. This asymmetric graph has a point symmetric outer shape (Fig. 3b). Therefore the graph in Figure 3a is asymmetric, but not completely asymmetric.}
  \end{figure}
  
  \newpage
  
  \section{\large{The smallest known completely asymmetric $(4;n)$-regular matchstick graphs for $4\leq n\leq11$}}
  
  \begin{figure}[!ht]
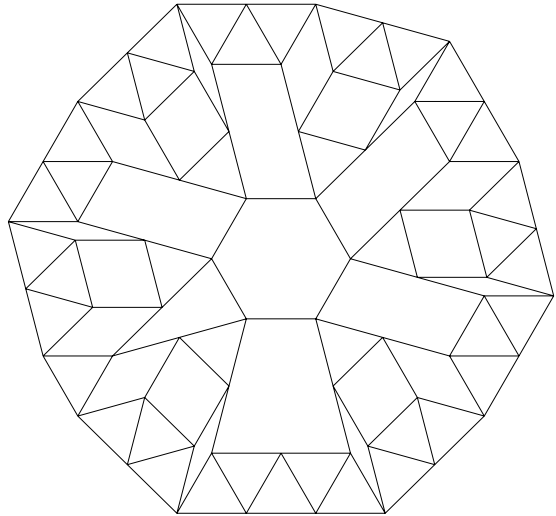

  	\centering
  	%
    \caption{4-regular matchstick graph with 66 vertices and 132 edges.}
  \end{figure}
  
  \begin{center}\end{center}
  
  \begin{figure}[!ht]
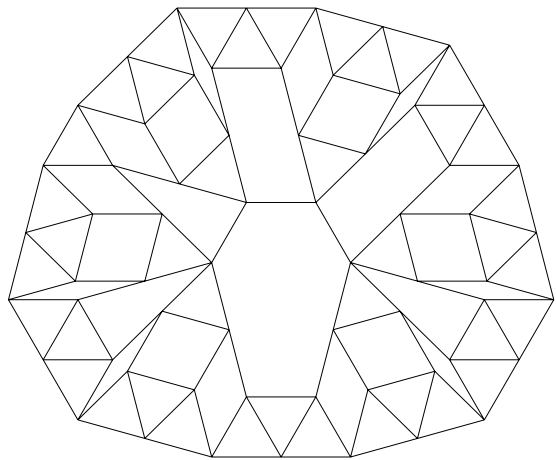

  	\centering
  	%
    \caption{$(4;5)$-regular matchstick graph with 62 vertices and 125 edges.}
  \end{figure}
  
  \newpage
  
  \begin{figure}[!ht]
  	\centering
  	%
    \caption{$(4;6)$-regular matchstick graph with 63 vertices and 128 edges - v1.}
  \end{figure}
  
  \begin{center}\end{center}
  
  \begin{figure}[!ht]
  	\centering
  	%
    \caption{$(4;6)$-regular matchstick graph with 63 vertices and 128 edges - v2.}
  \end{figure}
  
  \newpage
  
  \begin{figure}[!ht]
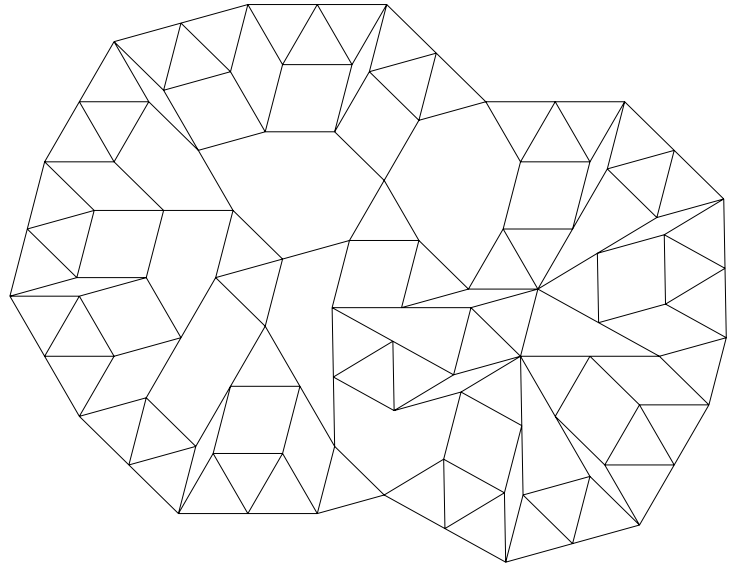

  	\centering
  	%
    \caption{$(4;7)$-regular matchstick graph with 93 vertices and 189 edges.}
  \end{figure}
  
  \begin{center}\end{center}
  
  \begin{figure}[!ht]
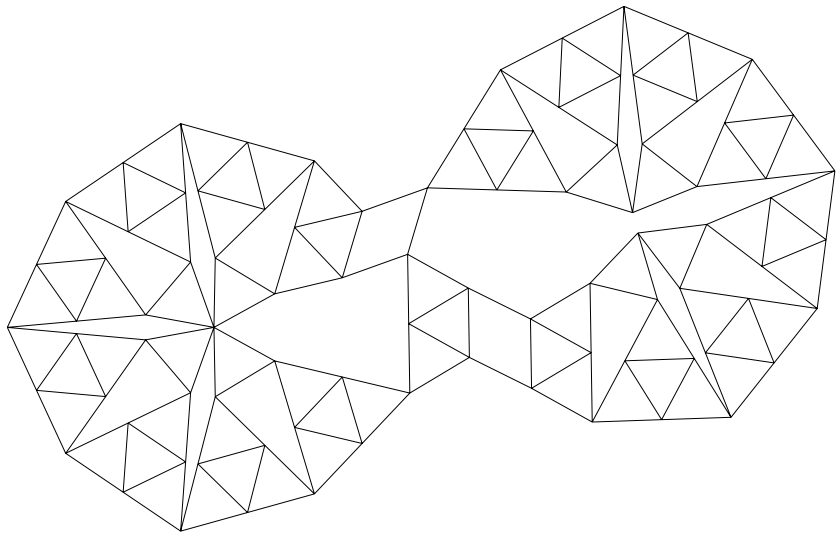

  	\centering
  	%
    \caption{$(4;8)$-regular matchstick graph with 87 vertices and 176 edges.}
  \end{figure}
  
  \newpage
  
  \begin{figure}[!ht]
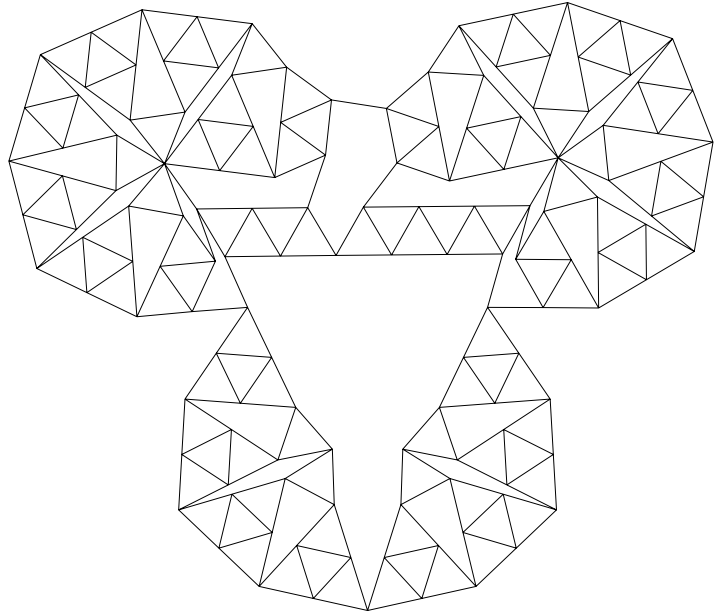

  	\centering
  	%
    \caption{$(4;9)$-regular matchstick graph with 136 vertices and 277 edges.}
  \end{figure}
  
  \begin{figure}[!ht]
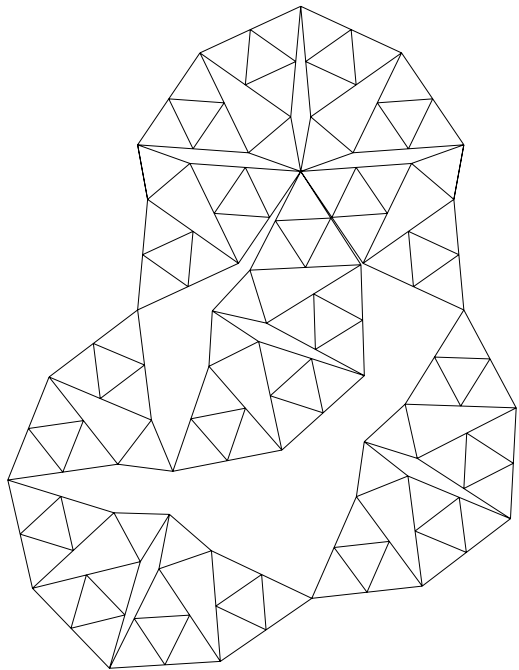

  	\centering
  	%
    \caption{$(4;10)$-regular matchstick graph with 114 vertices and 231 edges.}
  \end{figure}
  
  \newpage
  
  \begin{figure}[!ht]
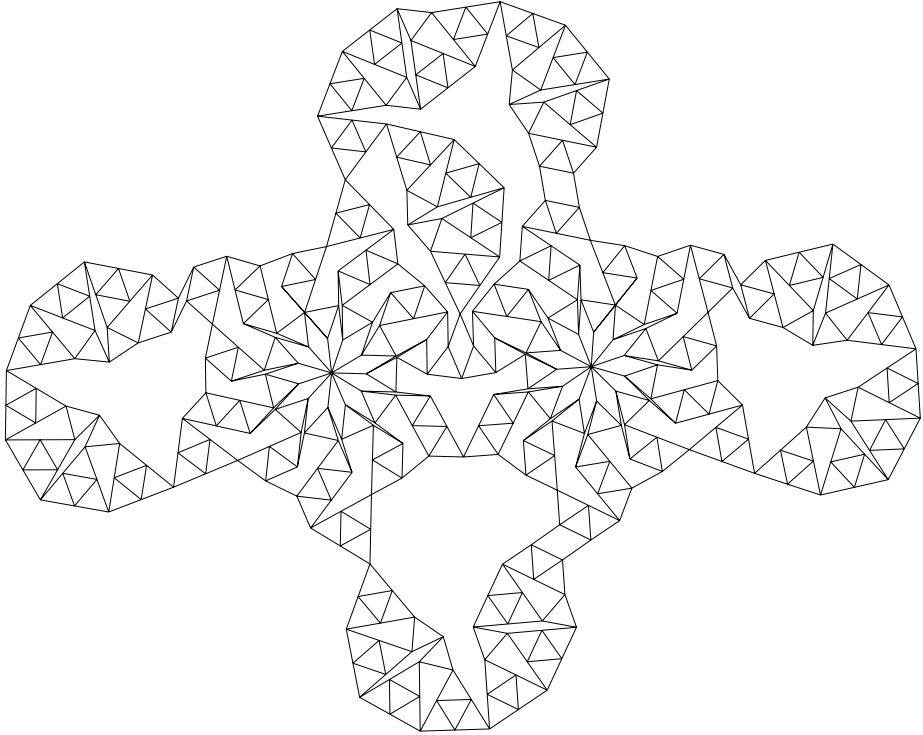

  	\centering
  	\caption{$(4;11)$-regular matchstick graph with 382 vertices and 771 edges.}
  	%
  \end{figure}
  \noindent 
  
  \newpage
  
  \noindent The interesting part of the graph in Figure 12 is the area around a vertex of degree 11, as the detail image in Figure 13 shows. This subgraph is rigid, but if removing only one edge the whole subgraph is flexible.
   
 \begin{figure}[!ht]
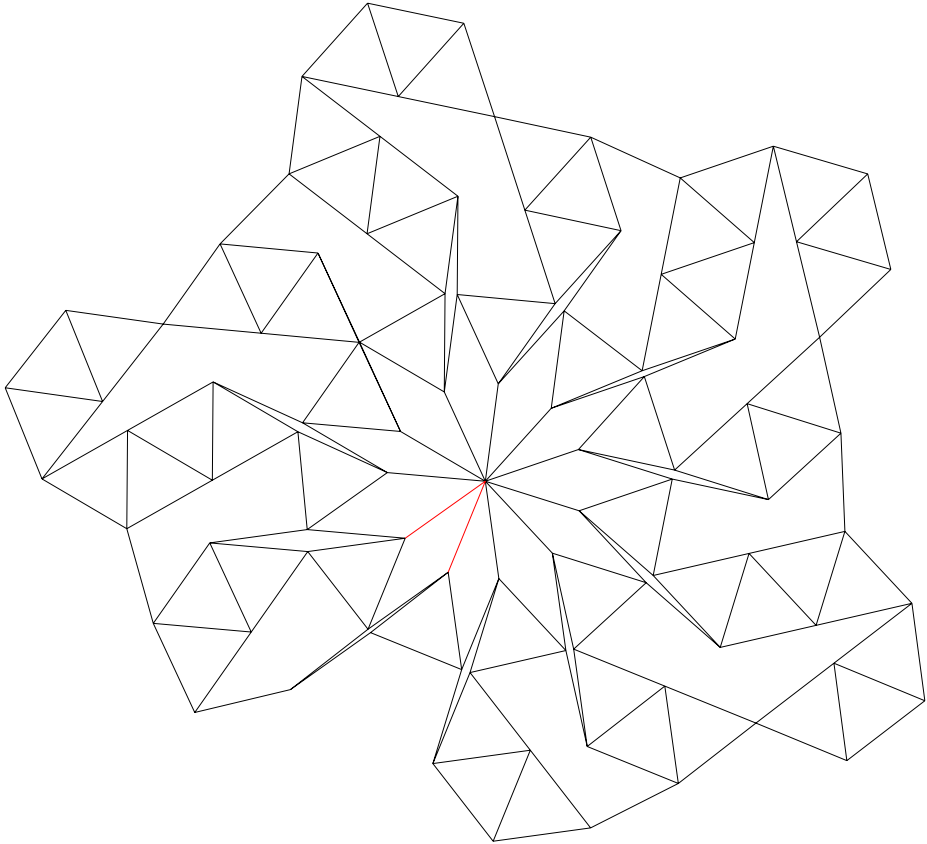

 	\centering
 	%
 	\caption{Detail around the right vertex of degree 11 in Figure 13.}
  \end{figure}
  
  \noindent This flexibility makes it possible to adjust the eleven angles around the centered vertex so that all edges have exactly one or two unit lengths. Beginning clockwise with the angle between the red edges, these degrees are\\
  \\
  32.362519660072210,\quad 40.49207000332465,\quad 25.382433534610843,\\
  34.890820876760450,\quad 32.21894760945070,\quad 34.514335947363630,\\
  29.108515978283318,\quad 36.31491131809427,\quad 29.550687898877964,\\
  35.065359484316880,\quad 30.09939768884507.
  
  \newpage
  
  \section{Remarks on the graphs}
  
  \noindent
  \textbf{Figure 4}: This $4$-regular matchstick graph with 132 edges was discovered on June 16, 2016 by M. Winkler. The graph has a triplet-kite-based geometry. There exist two other versions of this graph with a slightly different internal geometry, both with 134 edges.
  \\ \\
  \textbf{Figure 5}: This $(4;5)$-regular matchstick graph with 125 edges was discovered on June 25, 2016 by M. Winkler. The graph has a triplet-kite-based geometry and contains two triplet-kites. There exist seven other versions of this graph with a slightly different internal geometry, two with 126 edges, two with 127 edges, two with 128 edges, and one with 129 edges.
  \\ \\
  \textbf{Figure 6 and 7}: These $(4;6)$-regular matchstick graphs with 128 edges were discovered on June 25, 2016 by M. Winkler. Each graph has a triplet-kite-based geometry and contains two triplet-kites. Each version has a slightly different internal geometry.
  \\ \\
  \textbf{Figure 8}: This $(4;7)$-regular matchstick graph with 189 edges was discovered on June 16, 2016 by M. Winkler. The graph has a triplet-kite-based geometry and arises from a fusion of the graphs for $n=4$ and $n=5$.
  \\ \\
  \textbf{Figure 9}: This $(4;8)$-regular matchstick graph with 176 edges was discovered on August 12, 2016 by P. Dinkelacker. The graph has a kite-based geometry and contains two double-kites, two kites and two slightly modified kites.
  \\ \\
  \textbf{Figure 10}: This $(4;9)$-regular matchstick graph with 277 edges was discovered on August 11, 2016 by P. Dinkelacker. The graph has a kite-based geometry and contains four double-kites and four slightly modified kites. This graph is based on the currently smallest known symmetric equivalent with 273 edges.
  \\ \\
  \textbf{Figure 11}: This $(4;10)$-regular matchstick graph with 231 edges was discovered on March 17, 2016 by P. Dinkelacker. The graph has a kite-based geometry and consists only of kites. Three double kites, two reverse double kites and one kite. This version is also the currently smallest known one. It remains an interesting question whether a symmetric $(4;10)$-regular matchstick graph with 231 edges or less exists.
  \\ \\
  \textbf{Figure 12 and 13}: This $(4;11)$-regular matchstick graph with 771 edges was discovered on April 15, 2018 by M. Winkler, S. Vogel and P. Dinkelacker. This rigid graph is asymmetric and contains four double kites, one reverse double kite, four kites and four slightly modified kites. There exists a few asymmetric variations of this graph with 771 edges, because the clasps can be varied. But the current design requires the least place in the plane.
  
  \newpage
  
  \section{References}
  
  1. Heiko Harborth, \textit{Match Sticks in the Plane}, The Lighter Side of Mathematics. Proceedings of the Eug\`ene Strens Memorial Conference of Recreational Mathematics \& its History, Calgary, Canada, July 27 -- August 2, 1986 (Washington) (Richard K. Guy and Robert E. Woodrow, eds.), Spectrum Series, The Mathematical Association of America, 1994, pp. 281--288.
  \\ \\
  2. Sascha Kurz and Giuseppe Mazzuoccolo, \textit{3-regular matchstick graphs with given girth}, Geombinatorics Quarterly Volume 19, Issue 4, April 2010, pp. 156--175.\\
  \footnotesize(http://arxiv.org/pdf/1401.4360v1.pdf)\normalsize
  \\ \\
  3. Mike Winkler, Peter Dinkelacker, and Stefan Vogel, \textit{Streichholzgraphen 4-regul\"ar und 4/n-regul\"ar (n$>$4) und 2/5}, thread in a graph theory internet forum, used nicknames: P. Dinkelacker (haribo), M. Winkler (Slash),\\
  \footnotesize(https://tinyurl.com/ya3g6p7w)\normalsize
  \\ \\
  4. Stefan Vogel, \textit{Beweglichkeit eines Streichholzgraphen bestimmen}, July 2016.\\
  \footnotesize(https://tinyurl.com/yc8at6r7)\normalsize
  \\ \\
  5. Stefan Vogel, \textit{Matchstick Graphs Calculator (MGC)}, a software for the construction and calculation of matchstick graphs, 2016 -- 2017.\\
  \footnotesize(http://mikewinkler.co.nf/matchstick\_graphs\_calculator.htm)\normalsize
  \\ \\
  6. Mike Winkler, Peter Dinkelacker, and Stefan Vogel, \textit{New minimal (4; n)-regular matchstick graphs}, Geombinatorics Volume XXVII, Issue 1, July 2017, Pages 26--44, arXiv:1604.07134[math.MG].

\end{document}